\documentclass[12pt,a4article]{article}

\usepackage[latin1]{inputenc}
\usepackage{amsthm, amssymb, amscd,amsmath, amsfonts}
\usepackage{color}
\usepackage{graphicx}
\usepackage{hyperref}
\usepackage{authblk}
\usepackage{geometry}
\usepackage[T1]{fontenc}
 \usepackage[titletoc,toc,title]{appendix}

\usepackage{color,soul}
\usepackage{bigints}
\usepackage{fancyhdr}

\usepackage{hyperref}

\usepackage{amssymb}

\usepackage{mathrsfs,mathtools}
\usepackage[english]{babel}

\usepackage{amsmath,amsthm, bbm}

\usepackage{dsfont}
\usepackage{varioref}

\hypersetup{colorlinks=true, citecolor=blue, urlcolor=blue,linkcolor=blue}

\newtheorem{theorem}{Theorem}

\newtheorem{remark}[theorem]{Remark}

\def\R{\mathbb{R}}

\def\wh{\widehat}
\def\N{\mathbb{N}}

\def\E{\mathbb{E}}
\def\G{\mathscr{G} }
\def\to{\rightarrow}

\def\lto{\longrightarrow}

\def\bv{\big\vert}
\def\Bv{\Big\vert}
\def\1{\mathbbm{1}}
\def\e{\varepsilon}

\def\X{\mathfrak{X}}

\author{A. Lechiheb, I. Nourdin, G. Zheng and E. Haouala}

\title{Convergence of random oscillatory integrals in the presence of long-range dependence and application to homogenization}
\begin{document}

\maketitle

\begin{abstract}
This paper deals with the asymptotic behavior of random oscillatory integrals in the presence of long-range dependence. 
As a byproduct, we solve the corrector problem
in random homogenization of one-dimensional elliptic equations with highly oscillatory random coefficients displaying long-range dependence, by proving convergence to stochastic integrals with respect to Hermite processes.
\end{abstract}

 \section{Main results of the paper}

\subsection{Convergence of random oscillatory integrals}\label{S0}

In the present paper  one of our goals is to study, once properly normalized, the distributional convergence of some random oscillatory  integrals of the form
 \begin{align}\label{motive} 
 \int_0^1 \Phi\big[ g(x/\e) \big] h(x) \, dx  \,\, ,
 \end{align}
where
\begin{itemize}
\item  $h\in C\big([0,1]\big)$ is deterministic,
\item  $\big\{g(x)\}_{x\in\mathbb{R}_+}$ is a certain centred stationary Gaussian process exhibiting long-range correlation,
\item $\Phi\in L^2(\R,\nu)$ has Hermite rank $m\geq 1$ (with $\nu$ the standard Gaussian measure).
\end{itemize}

As we will see later, the main motivation of this study comes from the random corrector problem studied in \cite{GB12}.

Let us first introduce the
Gaussian process $\big\{ g(x) \big\}_{x\in\R_+}$  
we will deal with throughout all this paper. It is
constructed as follows:
  \begin{itemize}
 \item[1.] Let $m\in\N^\ast$ be fixed, let $H_0\in( 1 - \frac{1}{2m}, 1)$, and set $H = 1 + m(H_0 - 1)\in (1/2, 1)$;

  \item[2.] Fix  a slowly varying function $L: (0,+\infty)\to (0,+\infty)$  at $+\infty$, that is, consider a measurable and locally bounded function $L$ such that $L(\lambda x)/L(x)\to 1$ as $x\to+\infty$, for every $\lambda > 0$. Assume furthermore that $L$ is bounded away from $0$ and $+\infty$ on every compact subset of $(0,+\infty)$. (See \cite{BGT87} for more details on slowly varying functions.)

   \item[3.]  Let $e:\R\to \R$ be a square-integrable function
   such that
   \begin{itemize}
     \item [(3a)] $ \int_\R e(u)^2 \, du = 1 $,
     \item[(3b)] $|e(u)|\leq C u^{H_0-\frac32}L(u)$ for almost all $u>0$, for some absolute constant $C$,
     \item[(3c)] $e(u) \sim C_0 u^{H_0 -\frac{3}{2}} L(u)$, where $C_0 = \big( \int_0^\infty (u + u^2 )^{H_0-\frac{3}{2}} \, du \big)^{-1/2}$,
     \item[(3d)] their exist $0<\gamma < \min\big\{  H_0 - (1 - \frac{1}{2m} ) ,  1 - H_0    \big\}$  such that
      \begin{equation*}
    \int_{-\infty}^0 \vert e(u) e(xy+u)\vert \,du = o(x^{2H_0-2} L(x)^2) y^{2H_0 - 2 - 2\gamma} \end{equation*}
      as $x\to\infty$, uniformly in $y\in(0, t]$ for each given $t > 0$.
      \end{itemize}

     \item[4.] Finally, let  $W$ be a two-sided Brownian motion.
    
 \end{itemize}

Bearing all these ingredients in mind, we can now set, for $x\in\R_+$,
 \begin{equation}\label{g}
 g(x) := \int_{-\infty}^\infty e(x- \xi) dW_\xi \,\, .
 \end{equation}
 
 \begin{remark}
 {\rm 
    \begin{itemize}

    \item[(i)] Assumptions 3a and 4 ensure that $\big\{g(x)\big\}_{x\in\R_+}$ is a normalised centred Gaussian process.

     \item[(ii)]  Assumption 3b controls $\vert e(u) \vert$ for small $u$, while Assumption 3d ensures that the ``forward" contribution of $e(u)$ is ultimately negligible due to the following computation:
     \begin{align*}
     \E\big[  g(s) g(s+x)\big] & = \int_{-\infty}^\infty e(s- \xi) e(s+x -\xi)\, d\xi = \int_{-\infty}^\infty e(u) e(u+x)\, du \\
     & =  \int_{-\infty}^0 e(u) e(u+x)\, du  +  \int_{0}^\infty e(u) e(u+x)\, du \\
     & = o\big(x^{2H_0 - 2}L(x)^2\big) + x \int_0^\infty e(xu)e(xu +x)\, du \,\, .
     \end{align*}

\item[(iii)] Assumption 3c ensures that the process $\big\{ g(x)\big\}_{x\in\R_+}$ exhibits the following asymptotic behaviour:
      \begin{equation}
  R_g(x) := \E\big[ g(s) g(s+x) \big] \sim  x^{2H_0-2} L(x)^2 \quad\text{ as $x\to+\infty$},
   \end{equation}
see {\rm\cite[Equation (2.3)]{Taqqu79}}.
 \end{itemize}

}

\end{remark}

In section \ref{S2}, we will show that the random integral given by (\ref{motive}) exhibits the following asymptotic behavior as $\e\to 0$.

\begin{theorem}\label{atef}
Let $g$ be the centred stationary Gaussian process defined by (\ref{g}), and assume that $\Phi\in L^2(\R,\nu)$ has Hermite rank $m\geq 1$.
Then, for any   $h\in C\big([0,1]\big)$, the following convergence in law takes place
\begin{equation}\label{goal}
 M^\e_h : =\frac{1}{ \e d(1/\e) }\int_0^1 \Phi[g(x/\e)]  h(x)\,  dx \xrightarrow{\e\downarrow 0}   M^0_h := \frac{ V_m}{m!}  \int_0^1 h(x)\, dZ(x) \, ,
  \end{equation}
where $Z$ is the $m$th-Hermite process defined by (\ref{Z}) and $d(\cdot)$ is defined by
\begin{equation}\label{d(x)}
d(x) = \sqrt{\frac{m!}{H(2H-1)}} x^H L(x)^m.
\end{equation}

\end{theorem}

As we already anticipated, the fine analysis of the asymptotic behavior of (\ref{goal}) is motivated by the random corrector problem studied in \cite{GB12},
that we describe now.

\subsection{A motivating example}\label{S1}

Theorem \ref{atef} appears to be especially useful and relevant  in the study of  the following homogenization problem. Consider the following \emph{one-dimensional}
elliptic equation displaying random coefficients:
 \begin{align}\label{EE}
 \begin{cases}
 \displaystyle{ -\frac{d}{dx}\left( \,\, a(x/\e, \omega) \frac{d}{dx} u^\e(x, \omega) \, \right) } = f(x) \,\, , \quad x\in (0,1) \,\,, \quad\e > 0   \\
 u^\e(0, \omega) = 0 \,\,,\quad u^\e(1, \omega) = b\in\R.
 \end{cases}
 \end{align}
In (\ref{EE}),  the random potential $\{a(x)\}_{x\in\R_+}$ is assumed to be a uniformly bounded, positive\footnote{That is, there exists $r\in(0,1)$ such that  $r \leq a(x) \leq r^{-1}$ for every $(x, \omega)\in\R_+\times \Omega$.} stationary  stochastic process, whereas the data $f$ is continuous.  This model has received a lot of interests in the literature (see for instance \cite[page 13-14]{JKO94}).

Taking strong advantage of the fact that the ambient dimension is one, it is immediate to check that
the solution to  \eqref{EE} is given explicitly by
 \begin{equation}\label{ue}
 u^\e(x, \omega) = c^\e(\omega) \int_0^x \frac{1}{a(y/\e,\omega)} \, dy - \int_0^x \frac{F(y)}{a(y/\e, \omega)} \, dy,
 \end{equation}
 where
 $F(x) := \int_0^x f(y) \, dy$
 is the antiderivative of $f$ vanishing at zero,
 and where
 $$c^\e(\omega) : = \left( \, b + \int_0^1 \frac{F(y)}{a(y/\e, \omega)} \, dy \,\right) \left(  \int_0^1 \frac{1}{a(y/\e, \omega)} \, dy  \right)^{-1}.$$
Under suitable ergodic and stationary assumptions on $a$, the ergodic theorem applied to (\ref{ue}) implies that
$u^\e$ converges pointwise to $\bar{u}$ as $\e\to 0$, where
$$
\bar{u}(x) =  \frac{c^\ast x}{a^\ast} - \int_0^x \frac{F(y)}{a^\ast} \, dy,
$$
with $
 c^\ast : = b a^\ast + \int_0^1 F(y) \, dy
$
 and $$a^\ast: = \frac{1}{\E\big[  1/a(0)  \big] }.$$
 The above parameter $a^\ast$ is usually refered to as the \emph{effective diffusion coefficient} in the literature, see \emph{e.g.} \cite{PV81}.
It is also immediately checked that $\bar{u}$ is the unique solution to the
following \emph{deterministic} equation:
  \begin{align}\label{EEhomog}
 \begin{cases}
 \displaystyle{ -\frac{d}{dx}\left( \,\, a^\ast  \frac{d}{dx} \bar{u}(x) \, \right) } = f(x) \,\, , \quad x\in  (0,1) \\
 \bar{u}(0) = 0 \,\,,\quad \bar{u}(1) = b.
 \end{cases}
 \end{align}
Interested readers can refer to \cite{BG15} for a recent review on models involving more general elliptic equations.

In this work, we address the random corrector problem for (\ref{EE}) in presence of \emph{long-range} media, that is, we analyze the behaviour of the random fluctuations between $u^\e$ and $\bar{u}$ when the random potential $a$ is obtained by means of a long-range process (see below for the details). Taking advantage of the explicit expressions for both (\ref{EE}) and (\ref{EEhomog}), it is easy but crucial to observe that the random corrector $u^\e(x)-\bar{u}(x)$ can be fully expressed by means of random oscillatory  integrals of the form
\begin{equation}\label{integrales}
\int_0^1 \left[\,\, \frac{1}{a(y/\e)} - \frac{1}{a^\ast} \, \right] h(y) \, dy
\end{equation}
 for some  function $h$.   Thus, the random corrector problem for (\ref{EE}) reduces in a careful analysis of the asymptotic behaviour of random quantities of the form (\ref{integrales}) as $\e\to 0$. To this aim,
we need to give a precise description about the form of
the process $a$.

Let $\nu$ denote the standard Gaussian measure on $\R$. Every $\Phi\in L^2(\R, \nu)$ admits the following series expansion
 \begin{align}\label{chaotic}
  \Phi = \sum_{q=0}^\infty  \frac{V_q}{q!} H_q,  \quad\text{with $V_q:=  \int_\mathbb{R} \Phi(x)  H_q(x) \nu(dx)$,}
 \end{align}
 and where $H_q(x)= (-1)^q \exp(x^2/2)\frac{d^q}{dx^q} \exp(-x^2/2)$ denotes the $q$th Hermite polynomial. Recall that the integer $m_\Phi : = \inf\{q\geq 0:\,V_q\neq 0\}$ is called the {\it Hermite rank} of $\Phi$ (with the convention $\inf\emptyset = +\infty$). For any integer $m\geq 1$, we define $\mathscr{G}_m$ to the collection of all square-integrable functions (with respect to the standard Gaussian measure on $\R$) that have Hermite rank $m$.


Using Theorem \ref{atef} as main ingredient, we will  prove the following result about the asymptotic behaviour of the random corrector associated with (\ref{EE}).

\begin{theorem}\label{main}
Fix an integer $m\geq 1$ as well as two real numbers $H_0\in( 1 - \frac{1}{2m}, 1)$ and $b\in\R$,
and let $\{a(x)\}_{x\in\R_+}$ be a uniformly bounded, positive and stationary  stochastic process. Assume in addition that $q=\{q(x)\}_{x\in\R_+}$ given by 
\begin{equation}\label{a}
 q(x) =\dfrac{1}{a(x)} - \dfrac{1}{a^\ast},\quad \mbox{where }a^\ast: = 1/\E\big[  1/a(0)  \big],
\end{equation}
has the form 
\begin{equation}\label{q}
q(x) = \Phi\big( g(x) \big) \, ,
\end{equation}
 where $\Phi\in L^2(\R,\nu)$ belongs to $\mathscr{G}_m$
  and $\{ g(x) \}_{x\in\R_+}$ is the Gaussian process  given by (\ref{g}).
  Finally, let $f:[0,1]\to\R$ be continuous, and let us consider the solutions $u^\e$ and $\bar{u}$ of (\ref{EE}) and (\ref{EEhomog}) respectively.
 Then, for each $\e>0$,
the random corrector $u^\e-\bar{u}$ is a continuous process on $[0,1]$. Moreover,  we have the following convergence in law on $C([0,1])$ endowed with the supremum norm as $\e\to 0$:
     \begin{equation*}
     \left\{\frac{u^\e(x)-\bar{u}(x)}{ \e d(1/\e) } \right\}_{x\in[0,1]}\Longrightarrow  \,\, \left\{   \frac{V_m}{m!}  \int_\R F(x,y)\, dZ(y)\right\}_{x\in[0,1]} \,,
       \end{equation*}
where $d$ is given by (\ref{d(x)}),
\begin{align}
 F(x)&=\int_0^x f(y)dy,\quad c^\ast = a^\ast b + \int_0^1 F(y) \, dy,
\notag\\
 F(x,y) &=  \big[c^\ast  - F(y) \big]  \textbf{1}_{[0, x]}(y)  +  x \Big(   F(y) - \int_0^1 F(z) dz - a^\ast b  \Big)\textbf{1}_{[0, 1]}(y),\notag
 \end{align}
 and  $Z$ is the Hermite process of order $m$ and self-similar index $$H := 1+m(H_0 - 1)\in(1/2, 1) .$$
  (The definition of $Z$ is given in Theorem \ref{Taqqu79}.)
 \end{theorem}

Note that it is not difficult to construct a process $a$ satisfying all the assumptions of Theorem \ref{main}. 
Indeed, bearing in mind the notation of Theorem \ref{main}, we can write
\begin{equation}\label{ahah}
a(x)  = \left( q(x) + \frac{1}{a^\ast} \right)^{-1} =  \left( \Phi(g(x)) + \frac{1}{a^\ast} \right)^{-1} .
\end{equation}
Firstly, we note that since $g$ given by (\ref{g}) is stationary, clearly the same holds for $a$, whatever the expression of $\Phi$. Secondly, given any fixed $a^\ast > 0$, we can construct a bounded measurable function $\Phi\in\mathscr{G}_2$ with $\| \Phi \| _\infty \leq 1/(2a^\ast)$:
 \begin{itemize}
  \item[]  let $h_1, h_2$ be two bounded measurable functions, then it is clear that they belong to $L^2(\R, \nu)$ and they admit the series expansion
\begin{align*}
h_1 - \int_\mathbb{R} h_1 \, d\nu = \sum_{k=1}^\infty a_k H_k \, \qquad \text{and} \qquad h_2 - \int_\mathbb{R} h_2 \, d\nu = \sum_{k=1}^\infty b_k H_k \,\, ,
\end{align*}
where the coefficients $a_k, b_k$ are defined in the obvious manner. Therefore, the function 
 $$ \Psi : =  b_1 \left( h_1 - \int_\mathbb{R} h_1 \, d\nu \right) - a_1 \left( h_2 - \int_\mathbb{R} h_2 \, d\nu \right) $$ 
is bounded and belongs to $\mathscr{G}_2$.  
\end{itemize}
Then we pick $\Phi = \dfrac{\Psi}{2a^\ast\| \Psi \| _\infty }\in\G_2$. Therefore $a(x)$ defined by \eqref{ahah} satisfies 
\begin{align}\label{Good_news}0 < \frac{2a^\ast}{3} \leq  a(x) \leq 2 a^\ast \,\, .\end{align}
Inductively, one can construct a  bounded measurable $\Phi$ with Hermite rank $m\geq 3$ (by starting with two bounded functions in $\G_{m-1}$) such that  the process $\big\{ a(x) , x\in\mathbb{R} \big\}$ given in \eqref{ahah} satisfies \eqref{Good_news}.

Yet another possibility of constructing such a process $\big\{ a(x) , x\in\mathbb{R} \big\}$ is stated (more explicitly) as follows: fix $0<t_1<\ldots<t_m$, and consider the unique $(m+1)$-uple $(b_0,\ldots,b_m)$  satisfying
\begin{equation}\label{system}
\left\{
\begin{array}{lll}
\sum_{l=0}^{m} b_l\,e^{-kt_l}=0&\mbox{for all $k\in\{0,\ldots,m-1\}$},\\
\\
\sum_{l=0}^m b_l\,e^{-mt_l}=1 \,\, .
\end{array}
\right.
\end{equation}
(Existence and uniqueness of a solution to (\ref{system}) is a consequence of a Vandermonde determinant.)
Now, consider any measurable function $\psi$ satisfying 
\begin{equation}\label{bound}
0\leq \psi\leq \frac{1}{2a^*\sum_{l=0}^m|b_l|}.
\end{equation} 
Since $\psi$ belongs obviously to $L^2(\R,\nu)$, it may be expanded in Hermite polynomials as $\psi = \sum_{k=0}^\infty a_k H_k$. 
We assume moreover that $a_m\neq 0$. (Existence of $\psi$ satisfying both (\ref{bound}) and $a_m\neq 0$ is clear by a contradiction argument.)
Now, let $$\Phi=\sum_{l=0}^m b_l P_{t_l}\psi,$$
where $P_t\psi(x)=\int_\R \psi(e^{-t}x+\sqrt{1-e^{-2t}}y)\nu(dy)$ is the classical
Ornstein-Uhlenbeck semigroup.
Due to (\ref{system}), it is readily checked that the expansion of $\Phi$
is
$$
\Phi = a_m H_m +\sum_{k=m+1}^\infty \left\{ \sum_{l=0}^m b_l e^{-kt_l}\right\}a_kH_k,
$$
so that $\Phi\in\mathscr{G}_m$. Moreover,
$$
\|\Phi\|_\infty \leq \sum_{l=0}^m \vert b_l\vert \|P_{t_l}\psi\|_\infty
\leq \|\psi\|_\infty \sum_{l=0}^m \vert b_l\vert \leq \frac{1}{2a^*}
$$
and $a$ given by (\ref{ahah}) is positive and bounded.
So, existence of a process $a$ satisfying  all the assumptions of Theorem \ref{main} is shown.

 Theorem \ref{main} should be seen as an extension and unified approach of the main results of \cite{GB12}, and it contains these latter as particular cases. More precisely,
 the case where the Hermite  rank of $\Phi$
  is $m=1$ corresponds to \cite[Theorem 2.5]{BGMP08} and involves the fractional Brownian motion  in the limit, whereas the case where the Hermite  rank of $\Phi$
  is  $m=2$ corresponds to  \cite[Theorem 2.2]{GB12} and involves the Rosenblatt process in the limit.
  Also, in their last section (entitled {\it Conclusions and further discussion}), the authors of \cite{GB12} pointed out that ``it is natural to ask what would happen if the Hermite rank of $\Phi$ was greater than 2''.
  Our Theorem \ref{main} answers this question, by showing (as was guessed by the authors of \cite{GB12}) that, in the case $m\geq 3$, the limit
  takes the form of an integral with respect to the Hermite process of order $m$.
 Finally, we would like to emphasize that our Theorem \ref{main}, even in the cases $m=1$ and $m=2$, is a strict extension of the results of \cite{GB12}, as we allow the possibility to deal with a slowly varying function $L$.
That being said, our proof of Theorem \ref{main} is exclusively based on the ideas and results contained in
the seminal paper \cite{Taqqu79}  and follows the strategy developed in \cite{GB12}.  In higher dimension, it is usually very hard to study the corrector theory due to the lack of explicit form of the solution.  In a recent work \cite{MN15, MO16}, the authors considered the discretised version of the corrector problem in higher dimension and they were able to study the scaling limit to some Gaussian fields. For more details, we refer the interested readers to these two papers and the references therein.

\subsection{Plan of the paper}

The rest of the paper is organized as follows. In Section 2, we give some preliminary results, divided into several subsections. Section 3 contains the proof of Theorems \ref{main} and \ref{atef}.

 \section{Preliminary results}

 Throughout all this section, we let all the notation and assumptions of Sections \ref{S0} and \ref{S1} prevail.

\subsection{Asymptotic behaviour of the covariance function of $q$}\label{cov_aq} \qquad

For $x\in\R$, set $R_q(x) = \E\big[ q(0)q(x)\big]$.  Also, recall that $m$ is the Hermite rank of $\Phi$. Then, proceeding in similar lines as that of \cite[Lemma 2.1]{GB12}, one can show that
\begin{equation}\label{correlation}
\big\vert R_q(x) \big\vert = \big(  o(1) +    V_m^2/m! \big)     L(\vert x \vert)^{2m}  \vert x\vert^{-2(1-H)} \, ,
\end{equation}
as $\vert x \vert \to +\infty$. Here $o(1)$ means that the term converges to zero when $x\to \infty$.

The asymptotic relation \eqref{correlation} implies  the existence of some absolute  constant $C$ satisfying
\begin{align}\label{correlation_particular}
\bv R_q(x) \bv \leq C\, L( \vert x \vert )^{2m}  \vert x \vert^{-2(1-H)}
\end{align}
for any $x\neq 0$.

\subsection{Taqqu's theorem and convergence to Hermite process $Z$}
Recall $d(x)$ from (\ref{d(x)}). Its  main property  is that the variance of ${\displaystyle \dfrac{1}{d(x)}\int_0^x H_m(g(y))\, dy} $ is asymptotically equal to $1$ as $x\to +\infty$.

 The following result, due to Taqqu in 1979, is the key ingredient  in our proofs.


 \begin{theorem}({\rm \cite[Lemma 5.3]{Taqqu79})}\label{Taqqu79} \quad  Assume $\Phi\in\mathscr{G}_m$ and let $g$ be given by (\ref{g}). Then, as $T\to+\infty$,  the process
\begin{equation}\label{int-Taqqu}
 Y_T(x) = \frac{1}{d(T)} \int_0^{Tx}\Phi\big[ g(y) \big]  \, dy,\quad  x\in\R_+,
\end{equation}
 converges  to ${\displaystyle  \frac{V_m}{m!}  Z(x)}$ in the sense of finite-dimensional distributions, where  the $m$th-order Hermite process $Z$ with self-similar index $H =m(H_0 - 1) + 1$ is defined by:
\begin{align}
   Z(x)&\label{Z}\\
= K(m, H_0) & \left\{  \int_{-\infty}^\infty dB_{\xi_1}\int_{-\infty}^{\xi_1} dB_{\xi_2} \ldots \int_{-\infty}^{\xi_{m-1}} dB_{\xi_m}  \int_0^x \prod_{i=1}^m (s- \xi_i)^{H_0 - \frac{3}{2}}{\bf 1}_{(\xi_i < s)} \, ds   \right\} \,\, , \notag
\end{align}
where $$K(m, H_0): =  \sqrt{ \frac{m! H(2H-1) }{\Big( {\displaystyle\int_0^\infty (u + u^2 )^{H_0 - \frac{3}{2}} \, du}  \Big)^m}   }$$
is the normalising constant such that $\mathbb{E}\big[ Z(1)^2 \big] = 1$.  {\rm (See \cite[Equation (1.6)]{Taqqu79})}
\end{theorem}
Note that $Z(x)$ lives in the Wiener chaos of order $m$, which is non-Gaussian unless $m=1$ or $x=0$.

\subsection{Wiener integral with respect to $Z$}  \quad  Let $Z$ be given as above and let $\mathcal{E}$ be the set of elementary (deterministic) functions, that is, the set of functions  $h$ of the form
 $$ h(x) = \sum_{k=1}^\ell a_k \textbf{1}_{(t_k, t_{k+1}] }(x) $$
with $\ell\in\N^\ast$, $a_k\in\R$, $t_k < t_{k+1}$.  For such $h$, we define the Wiener integral with respect to $Z$ in the usual way, as a linear functional over $\mathcal{E}$:
 $$\int_\R h(x) \, dZ(x)  =  \sum_{k=1}^\ell a_k\Big[ Z(t_{k+1}) - Z(t_k)  \Big] \,\, .  $$
 One can verify easily that this definition is independent of choices of representation for elementary functions.   Now we introduce the space of (deterministic) integrands for this Wiener integral:
        \begin{equation}
  \Lambda^H = \left\{ \,\, f:  \R \lto\R  \, \Big\vert\,  \quad\int_\R \int_\R f(u)f(v) \vert u-v \vert^{2H-2} \, du\, dv  <  + \infty  \,\,  \right\},
            \end{equation}
equipped with the norm
\begin{equation}
  \| f \|^2_{\Lambda^H} = H(2H-1)\int_\R\int_\R   f(u)f(v) \vert u-v\vert^{2H-2}\, du \, dv \, .
\end{equation}
When $h\in\mathcal{E}$, it is straightforward to check the following isometry property:
           $$ \E\left[\left(\int_\R h(x)dZ(x)\right)^2 \right] = \|h\|^2_{\Lambda^H}. $$
As a consequence, one can define the Wiener integral $\int_\R f(x)dZ(x)$  for any $f\in \Lambda^H$, by a usual approximation procedure.

It is by now well known (thanks to  \cite{PT00}) that  $\big( \Lambda^H, \| \cdot \| _{\Lambda^H} \big)$ is a Hilbert space that contains distributions in the sense of Schwartz. To overcome this problem, we shall restrict ourselves to the proper subspace
 $$\vert {\Lambda^H} \vert =  \left\{ \,\, f:\R\to\R  \, \Big\vert\,  \quad\int_\R\int_\R \vert f(u)f(v) \vert  \vert u-v \vert^{2H-2} \, du\, dv  <  + \infty  \,\,  \right\}  $$
equipped with the norm
$$   \|f\|^2_{\vert  \Lambda^H  \vert} = H(2H-1) \int_\R\int_\R \vert f(u)f(v)\vert \vert u-v\vert^{2H-2}\, du \, dv \, .$$
We then have (see \cite[Proposition 4.2]{PT00})
\begin{equation}
 L^{1} \big(\R \big)\cap L^{2} \big(\R \big)\subset L^{1/H} \big(\R \big)\subset \vert \Lambda^H \vert \subset  \Lambda^H .\label{inclusion}
\end{equation}
Moreover, $\big( \vert \Lambda^H \vert,   \| \cdot \|_{\vert \Lambda^H \vert} \big)$ is  a Banach space, in which the set $\mathcal{E}$  is dense. So for $h\in  \vert \Lambda^H \vert$, we can define
\begin{equation}\label{cv}
\int_\R h(x) \, dZ(x) = \lim_{n\to+\infty}\int_\R h_n(x) \, dZ(x) \,,
\end{equation}
where $(h_n)$ is any sequence of $\mathcal{E}$ converging to $h$ in  $\big(  \vert \Lambda^H \vert,   \|\cdot \|_{ \vert \Lambda^H \vert} \big)$; the convergence in (\ref{cv}) takes place in $L^2(\Omega\big)$.

  For a detailed account of this integration theory, one can refer to \cite{MT07,PT00}.

\subsection{Some facts about slowly varying functions} \quad  Let $L: (0,+\infty) \to (0, +\infty)$ be a slowly varying function at $+\infty$ and $\alpha > 0$. It is well known (see \cite[Proposition 1.3.6(v)]{BGT87}) that
\begin{align*}
x^\alpha L(x) \to  + \infty \,  \quad\text{and} \quad x^{-\alpha} L(x) \to 0 \,\,,
\end{align*}
as $x\to +\infty$. In particular, one can deduce  that
\begin{align} \label{SL_fact1}
\lim_{\e\downarrow 0} \e^{1-H} L(1/\e)^m = 0 \,\, .
\end{align}

 The following result is known as Potter's Theorem (see \cite[Theorem 1.5.6(ii)]{BGT87}).

\begin{theorem}\label{Harry_Potter} Let $L: (0,+\infty) \to (0, +\infty)$ be a slowly varying function at $+\infty$ such that  it  is bounded away from $0$ and $+\infty$ on every compact subset of $(0,+\infty)$.  Then for any $\delta > 0$, there exists some constant $C = C(\delta)$ such that
 $$\frac{L(y)}{L(x)} \leq C \max\Big\{ \, (x/y)^\delta \,, \, (y/x)^{\delta} \Big\}  $$
for any $x, y\in(0,+\infty)$.

\end{theorem}

\section{Proofs of main results}

\subsection{Proof of Theorem \ref{atef}}\label{S2}

First recall that a typical function $h$ in $\mathcal{E}$ has the form
\begin{equation*}
  h(x)=\sum_{\ell=1}^n   a_\ell \textbf{1}_{(t_\ell,t_{\ell+1}]}(x),\quad t_\ell < t_{\ell+1},\quad a_\ell\in\R,\quad \ell=1,...,n\, .
\end{equation*}
For such a simple function $h$, we deduce from Taqqu's Theorem \ref{Taqqu79} that
\begin{align*}
  M^\e_h & = \frac{1}{\e d(1/\e) } \int_\R q(x/\e)\sum_{\ell=1}^n a_\ell \textbf{1}_{(t_\ell,t_{\ell+1}]}(x) \, dx    \\
  & = \sum_{\ell=1}^n a_\ell  \,\frac{1}{d(1/\e)} \left( \,\,   \int_0^{t_{\ell+1}/\e}  \Phi(g(x)) \,  dx  -  \int_0^{t_\ell  /\e }  \Phi(g(x))  \, dx \,\, \right)       \\
   &\xrightarrow{ \e \to 0 } \frac{V_m}{m!}\sum_{\ell =1}^n   a_\ell     \big[    Z(t_{\ell+1}) - Z(t_\ell)  \big]   =  \frac{V_m}{m!} \int_\R h(x)\,  dZ(x)\, .
\end{align*}
This proves  (\ref{goal}) for simple functions $h\in\mathcal{E}$.

Let us now  consider $h\in C\big( [0,1] \big)$. It is easy to see that there exists a sequence $(h_n)\subset\mathcal{E}$ such that
                             $$ \lim_{n\to+\infty} \big\| h_n - h \big\| _\infty = 0 \,\, .$$
Let us fix a number $\zeta\in (0,1)$ and  show the  convergence in $L^2(\Omega)$ of $M^\e_{h_n}$, uniformly in $\e\in(0, \zeta)$. First, one can write
\begin{align}
 &\qquad \sup_{\e\in(0,\zeta)} \E\Big[ \, \vert M^\e_{h_n}-M^\e_{h}\vert^2 \, \Big]  \nonumber \\
 &=\sup_{\e\in(0,\zeta)} \frac{1}{\e^2 d(1/\e)^2 } \,\E\left[ \,\, \left\vert \int_0^1 q(x/\e) \big[ h_n(x)-h(x) \big]  \, dx \, \right\vert^2   \,\,   \right] \nonumber \\
   & \leq \big\| h_n - h \big\| _\infty^2 \times \sup_{\e\in(0,\zeta)}  \frac{1}{\e^2 d(1/\e)^2}    \int_{\R^2\setminus D} \Big\vert R_q\left( \frac{y-x}{\e}\right) \Big\vert   \, dx\, dy \,\, ,  \nonumber
\end{align}
where $D = \big\{ (x,y)\in [0,1]^2 : x =y \big\}$ is a negligible subset of $\R^2$.  By \eqref{correlation_particular},
     $$\bv R_q\big( \frac{y-x}{\e}\big) \bv \leq \text{Cst}\, L\left(\left| \frac{y-x}{\e} \right|\right)^{2m}  \left| \frac{y -x}{\e} \right|^{-2(1-H)} \,, \quad  \forall (x,y)\in \R^2\setminus D\,.$$
Secondly, with $\beta > 0$ small enough such that  $2m\beta + 2(1-H)\in(0,1)$, we have
\begin{align}
& \qquad \sup_{\e\in(0,\zeta)} \frac{1}{\X(\e)^2} \int_{[0,1]^2\setminus D} \Big\vert R_q\left( \frac{y-x}{\e}\right) \Big\vert   \, dx\, dy  \notag   \\
&\leq  \text{Cst} \sup_{\e\in(0,\zeta)}  \int_{[0,1]^2\setminus D}   \left\{ \,\, \frac{ L\big( \vert (x-y)/\e \vert \big) }{ L(1/\e)} \, \right\} ^{2m}  \vert x-y\vert^{-2(1-H)}   \, dx\, dy\notag \\
&\leq \text{Cst}   \int_{[0,1]^2\setminus D}   \vert x-y\vert^{-2m\beta-2(1-H)}   \, dx\, dy \label{Potter} \\
& < +\infty  \,\, ,  \notag
      \end{align}
where \eqref{Potter} follows  from Theorem \ref{Harry_Potter}.   It is now clear that, indeed,
\begin{equation}\label{unifbound}
  \lim_{n\to+\infty} \sup_{\e\in(0,\zeta)}\E \big[ \vert M^\e_{h_n}-M^\e_{h}\vert^2\big] = 0 \, .
\end{equation}
To conclude, let $d(\cdot,\cdot)$ denote any distance metrizing the convergence in distribution between real-valued random variables
(for instance, the Fortet-Mourier distance).
For $h\in C([0,1])$ and $(h_n)\subset\mathcal{E}$ converging to $h$, one can write, for any $\e>0$ and $n\in\mathbb{N}$:
$$
d(M^\e_{h},M^0_h)\leq d(M^\e_{h},M^\e_{h_n})+
d(M^\e_{h_n},M^0_{h_n})+ d(M^0_{h_n},M^0_h).
$$
Fix $\eta>0$. By (\ref{unifbound}), one can choose $n$ big enough so that, for
any $\e\in(0, \zeta)$, both $d(M^\e_{h},M^\e_{h_n})$ and $d(M^0_{h_n},M^0_h)$ are less than $\eta/3$.
It remains to choose $\e>0$ small enough so that
$d(M^\e_{h_n},M^0_{h_n})$ is less that $\eta/3$ (by
 (\ref{goal}) for the simple function $h_n\in\mathcal{E}$), to conclude that
 (\ref{goal}) holds true for any continuous function $h$.

\begin{remark}\label{rem}
{\rm
Clearly, the above result still holds true for any function $h$ that is continuous except at finitely many points. Note also that the function $\Phi\in\mathscr{G}_m$ is not necessarily bounded in Theorem \ref{atef}.
}
\end{remark}

\subsection{Proof of Theorem \ref{main}} \quad  The proof is divided into five steps.  We write $\X(\e) = \e d(1/\e) =\sqrt{  \frac{m!}{H(2H-1)} }  \e^{1-H}L(1/\e)^m$.

\bigskip

{\bf (a) Preparation}.  \quad Following \cite{GB12}, more precisely identities (5.1) and (5.19) therein,
we first rewrite the rescaled corrector as follows:
\begin{align}\label{deomp_0}
  \frac{u^\e(x)-\bar{u}(x)}{ \X(\e) } = \mathcal{U}^\e(x) + \underbrace{    \frac{1}{\X(\e)}     r^\e(x) +   \frac{1}{\X(\e)}  \rho^\e\frac{x}{a^\ast} }_{=:\mathcal{R}^\e(x)} \,\, ,
\end{align}
where
\begin{align*}
\mathcal{U}^\e(x) &= \frac{1}{\X(\e)} \int_\R F(x, y) q(y/\e)  \, dy \, ,\\
r^\e(x) &= (c^\e-c^\ast)\int_0^x q(y/\e) \, dy \, ,
\end{align*}
and
\begin{align*}
  \rho^\e & := \frac{a^\ast}{\int_0^1 a(y/\e)^{-1} \, dy}   \bigg[ \left( \, a^\ast b+\int_0^1F(y)dy \right)\left(\int_0^1 q(y/\e)  \, dy\right)^2 \\
  & \qquad\qquad\qquad\qquad - \int_0^1F(y) q(y/\e)  \, dy  \times \int_0^1 q(y/\e)  \, dy  \bigg].
\end{align*}

Now, let us first show the weak convergence of $\mathcal{U}^\e$ to $\mathcal{U}$ in $C([0,1])$ and then prove that $\mathcal{R}^\e$ is a remainder.
In order to prove the first claim, we start by  establishing the f.d.d. convergence and then we prove the tightness.

\bigskip

{\bf (b) Convergence of finite dimensional distributions of $\mathcal{U}^\e$.}  \quad For $x_1,  \ldots, x_n\in\R$ and $\lambda_1, \ldots, \lambda_n\in\R$ ($n\geq 1$), we have
   \begin{align*}
\sum_{k=1}^n \lambda_k  \,\mathcal{U}^\e(x_k)  =   \frac{1}{\X(\e)} \int_\R  \sum_{k=1}^n \lambda_k \, F(x_k, y) q(y/\e) \, dy.
   \end{align*}
Note that the function $\sum_{k=1}^n \lambda_k \, F(x_k, \cdot)$ have at most finitely many discontinuities. Thus, Theorem \ref{atef} and Remark \ref{rem} imply that
    $\sum_{k=1}^n \lambda_k \, \mathcal{U}^\e(x_k)$ converges in distribution to $\sum_{k=1}^n \lambda_k \,\mathcal{U}(x_k)$,
yielding the desired convergence of finite dimensional distributions.

\bigskip

{\bf (c) Tightness of $\mathcal{U}^\e$}. \quad We check Kolmogorov's criterion (\cite[Corollary 16.9]{OK02}).  First observe  that $\mathcal{U}^\e(0) = 0$. Now, fix  $0 \leq u < v \leq 1$, and
set $F_1(y) = c^\ast - F(y)$ and $F_2(y) = F(y) - \int_0^1 F(t)\, dt - a^\ast  b$, so that $F(x,y)=F_1(y)\textbf{1}_{[0,x]}(y)+xF_2(y)\textbf{1}_{[0,1]}(y)$.
Then
\begin{align}
&\qquad  \E \big( \vert \mathcal{U}^\e(u) - \mathcal{U}^\e(v) \vert^2 \big) \notag \\
& = \E \bigg[ \,\, \frac{1}{\X(\e)^2} \bigg\vert \int_0^1 \textbf{1}_{(u, v]}(y) q(y/\e)  F_1(y) \, dy + (v-u)  \int_0^1 q(y/\e)   F_2(y)\, dy \, \bigg\vert^2 \, \bigg]  \notag\\
& \leq \frac{2}{\X(\e)^2} \E \bigg[ \bigg\vert \int_0^1 \textbf{1}_{(u, v]}(y) q(y/\e) F_1(y) \, dy  \bigg\vert^2  +    \bigg\vert(v-u) \int_0^1 q(y/\e)   F_2(y)\, dy \,\, \bigg\vert^2 \bigg]  \notag\\
& \leq \frac{2}{\X(\e)^2} \int_u^v \int_u^v F_1(x) F_1(y)  R_q\Big( \frac{y-x}{\e} \Big) \, dx\, dy  \label{trouble1} \\
&  \qquad\qquad\qquad +   \frac{2(v-u)^2}{\X(\e)^2} \int_0^1\int_0^1 F_2(x) F_2(y)  R_q\Big( \frac{y-x}{\e} \Big) \, dx \, dy \,\, . \label{trouble2}
\end{align}
Note that $F_2$ is bounded on $[0,1]$. Therefore, as far as (\ref{trouble2}) is concerned, one can write, using Potter's Theorem as in the proof of Theorem \ref{atef},
\begin{align*}
 \sup_{\e\in(0,\zeta)} \left\vert \frac{(v-u)^2}{\X(\e)^2} \int_0^1\int_0^1 F_2(x) F_2(y) R_q\Big( \frac{y-x}{\e} \Big) \, dx \, dy\right\vert   \leq  \text{Cst} (v-u)^2 \, .
\end{align*}
Now, let us consider the term in \eqref{trouble1}. Similarly,
         \begin{align}
& \qquad   \sup_{\e\in(0,\zeta)}\frac{1}{ \X(\e)^2 }  \left\vert \int_u^v \int_u^v F_1(x) F_1(y) R_q\Big( \frac{y-x}{\e} \Big) \, dx\, dy\right\vert   \notag   \\
&\leq  \text{Cst} \sup_{\e\in(0,\zeta)}  \frac{1}{\X(\e)^2 }   \int_u^v \int_u^v \bv R_q( \frac{y-x}{\e} ) \bv\, dx\, dy  \quad\text{(since $F_1$ is  bounded)} \notag\\
&\leq  \text{Cst} \sup_{\e\in(0,\zeta)}  \frac{1}{ L(1/\e)^{2m}  }   \int_u^v \int_u^v L\big( \vert y-x\vert/\e \big)^{2m} \frac{dx\, dy }{\vert y - x\vert^{2(1-H)}}  \notag    \\
& \leq  \text{Cst}  \int_u^v \int_u^v \vert y - x\vert^{-2(1-H) - 2m\beta}\, dy\, dx  \quad\text{(similarly as in \eqref{Potter})}  \notag\\
& = \text{Cst} (v-u)^{2 - 2m\beta - 2(1-H)} \, .  \label{similarly_as_in}
\end{align}
Since $2-2m(1-H_0) - 2m\beta>1$, this proves the tightness  of $(\mathcal{U}^\e)_{\e}$ by means of the usual Kolmogorov's criterion.

\bigskip

{\bf (d) Control on the remainder term $\mathcal{R}^\e$ in \eqref{deomp_0}}.  \quad We shall prove that the process $\mathcal{R}^\e$ converges in probability to zero in $C([0,1])$. First we claim that if $G\in C([0,1])$, then there exists some constant $C = C(G)$ such that
 \begin{align} \label{claim0} \sup_{x\in[0,1]}  \E\left[\left(  \int_0^x q(y/\e) G(y) \, dy \right)^2\right] \leq  C\, \X(\e)^2 \,\, . \end{align}
Indeed, the same argument we used for bounding \eqref{trouble2} works here as well:
\begin{align*}
 &\qquad \sup_{x\in[0,1]}  \E\left[\left( \int_0^x q(y/\e) G(y) \, dy \right)^2\right]  \\
 & \leq \| G \| _\infty^2  \int_{[0,1]^2} \bv R_q( \vert y - z \vert/\e) \bv \, dy\, dz  \\
 & \leq   \| G \| _\infty^2   \X(\e)^2 \left( \sup_{\e\in(0,\zeta)} \frac{1}{\X(\e)^2} \int_{[0,1]^2} \bv R_q( \vert y - z \vert/\e) \bv \, dy\, dz \right) \\
 & \leq \text{Cst}  \X(\e)^2 \,\, ,
\end{align*}
where the last inequality follows from \eqref{Potter}.

Now, let us consider $\mathcal{R}^\e$:

           (i) Due to the explicit expression of $\rho^\e$,  it follows from \eqref{claim0}, the fact that $a$ is bounded from below and Cauchy-Schwarz inequalities that
 \begin{align*}
 &\qquad \E\big[ \vert \rho^\e \vert \big] \\
  &\leq   \text{Cst} \left\{ \, \Big\| \int_0^1 q(y/\e) \, dy \Big\| _{L^2(\Omega)}^2 +   \Big\| \int_0^1 F(y)q(y/\e) \, dy \Big\| _{L^2(\Omega)}  \Big\| \int_0^1 q(y/\e) \, dy \Big\| _{L^2(\Omega)}  \,\right\} \\
  &\leq \text{Cst}\, \X(\e)^2 \,\, .
                   \end{align*}
   (ii)  Observe that
  \begin{eqnarray*}
  c^\e - c^\ast &=& a^\ast  \int_0^1\big( F(y) - \int_0^1 F(t) \,dt - ba^\ast \big) q(y/\e) \, dy + \rho^\e \\
  &=:&  \int_0^1\wh{F}(y) q(y/\e) \, dy + \rho^\e .
  \end{eqnarray*}
  Then
  \begin{align*}
   &\qquad \sup_{x\in[0,1]} \E\big[ \vert r^\e(x) \vert \big] =  \sup_{x\in[0,1]} \E\left[  \bv (c^\e - c^\ast)\int_0^x q(y/\e)\, dy  \bv \right] \\
  &\leq  \sup_{x\in[0,1]}\E \left[ \,  \Bv \int_0^1 \wh{F}(y) q(y/\e) \, dy  \int_0^x q(y/\e)\, dy \Bv  \,\right]  + \text{Cst}\, \E\big[ \vert \rho^\e\vert \big]  \leq  \text{Cst}\, \X(\e)^2  \,\, .
  \end{align*}
Therefore, as $\e\to 0$ we have
    $$\sup_{x\in[0,1]} \E\Big[ \bv \mathcal{R}^\e(x) \bv \Big] \leq \text{Cst} \X(\e) \lto 0 \,\, . \quad \text{(by \eqref{SL_fact1})} $$
In particular,  $\big\{\mathcal{R}^\e(x), x\in[0,1] \big\}$ converges to zero in the sense of finite-dimensional distributions.
Now, let us check the tightness of $\big(\mathcal{R}^\e\big)_\e$. Note that  $\mathcal{R}^\e(0)=0$ and that, for $0\leq u < v\leq 1$,
\begin{align*}
&\qquad \big\|  \mathcal{R}^\e(u) - \mathcal{R}^\e(v) \big\| _{L^2(\Omega)}^2 \\
&\leq \frac{2}{\X(\e)^2}  \left\{ \,\,  \big\| r^\e(u) - r^\e(v) \big\| ^2_{L^2(\Omega)} + \frac{2(u-v)^2}{\vert a^\ast\vert^2} \E\big[ \vert \rho^\e \vert^2 \big] \,\, \right\}  \\
& \leq   \frac{2}{\X(\e)^2}    \big\| r^\e(u) - r^\e(v) \big\| ^2_{L^2(\Omega)} \\
&\quad+\text{Cst} \frac{(u-v)^2}{\X(\e)^2} \E\big[ \vert \rho^\e \vert \big]\quad\mbox{(since $\rho^\e$ is uniformly bounded)}\\
&\leq   \frac{2}{\X(\e)^2}  \big\| r^\e(u) - r^\e(v) \big\| ^2_{L^2(\Omega)} + \text{Cst} (u-v)^2 \qquad \text{(by point (i) above)} \\
&\leq   \text{Cst}   \frac{1}{\X(\e)^2}  \int_{[u,v]^2} \bv R\big( (y-z)/\e\big) \bv \, dy\, dz  \\
&\quad + \text{Cst} (u-v)^2
\quad\text{(since $c^\e - c^\ast$ is uniformly bounded)}\\
&\leq   \text{Cst} (v-u)^{2 - 2(1-H) - 2m\beta} +  \text{Cst} (v-u)^2 \,,
\end{align*}
where the last inequality follows from the same arguments as in \eqref{similarly_as_in}.  Therefore, $\mathcal{R}^\e$ converges in distribution to $0$,  as $\e\downarrow0$, so it converges in probability to $0$.

\bigskip

{\bf (e) Conclusion}.
Combining the results of ({\bf a}) to ({\bf d}),
the proof of Theorem \ref{main} is  concluded by evoking Slutsky lemma.\qedhere

\end{document}